\documentclass[11pt]{article}
\usepackage{epsfig}

\def\be{\begin{equation}}
\def\ee{\end{equation}}
\def\bea{\begin{eqnarray}}
\def\eea{\end{eqnarray}}
\def\bes{\begin{eqnarray*}}
\def\ees{\end{eqnarray*}}

\def\nn{\nonumber}
\def\lb{\label}
\def\bs{\setminus}

\def\vs{{\varsigma}}

\def\R{{\bf R}}
\def\C{{\bf C}}
\def\Z{{\bf Z}}

\def\N{{\bf N}}
\def\U{{\bf U}}

\def\Q{{\bf Q}}

\def\aa{{\alpha}}
\def\bb{{\beta}}
\def\ga{{\gamma}}

\def\th{{\theta}}

\def\om{{\omega}}
\def\Om{{\Omega}}

\def\lm{{\lambda}}
\def\Lm{{\Lambda}}

\def\sg{{\sigma}}
\def\dm{{\diamond}}

\def\vf{{\varphi}}

\def\<{{\langle}}
\def\>{{\rangle}}

\def\P{{\cal P}}

\def\Sp{{\rm Sp}}

\def\ol{\overline}

\def\hb{\vrule height0.18cm width0.14cm $\,$}

\title{Two elliptic closed geodesics on positively curved Finsler spheres}

\author{Huagui Duan \thanks{Partially supported by NNSF (No.11131004, 11471169), LPMC of MOE of China and
Nankai University. E-mail: duanhg@nankai.edu.cn}\\\\
School of Mathematical Sciences and LPMC\\
Nankai University\\
Tianjin 300071, The People's Republic of China\\}

\begin{document}
\date{ }
\maketitle

\begin{abstract}
{\it In this paper, we prove that for every Finsler $n$-dimensional sphere $(S^{n},F)$ with reversibility $\lm$ and flag curvature $K$ satisfying $\left(\frac{\lm}{1+\lm}\right)^2<K\le 1$, either there exist infinitely many closed geodesics, or there exist at least two elliptic closed geodesics and each linearized Poincar\'{e} map has at least one eigenvalue of the form $e^{\sqrt{-1}\th}$ with $\th$ being an irrational multiple of $\pi$.}
\end{abstract}



\renewcommand{\theequation}{\thesection.\arabic{equation}}
\renewcommand{\thefigure}{\thesection.\arabic{figure}}

\setcounter{figure}{0}
\setcounter{equation}{0}
\section{Introduction and main result}

A closed curve on a Finsler manifold is a closed geodesic if it is
locally the shortest path connecting any two nearby points on this
curve. As usual, on any Finsler manifold
$(M, F)$, a closed geodesic $c:S^1=\R/\Z\to M$ is {\it prime}
if it is not a multiple covering (i.e., iteration) of any other
closed geodesics. Here the $m$-th iteration $c^m$ of $c$ is defined
by $c^m(t)=c(mt)$. The inverse curve $c^{-1}$ of $c$ is defined by
$c^{-1}(t)=c(1-t)$ for $t\in \R$.  Note that unlike Riemannian manifold,
the inverse curve $c^{-1}$ of a closed geodesic $c$
on a irreversible Finsler manifold need not be a geodesic.
We call two prime closed geodesics
$c$ and $d$ {\it distinct} if there is no $\th\in (0,1)$ such that
$c(t)=d(t+\th)$ for all $t\in\R$.
On a reversible Finsler (or Riemannian) manifold, two closed geodesics
$c$ and $d$ are called { \it geometrically distinct} if $
c(S^1)\neq d(S^1)$, i.e., their image sets in $M$ are distinct.
We shall omit the word {\it distinct} when we talk about more than one prime closed geodesic.

For a closed geodesic $c$ on $n$-dimensional manifold $(M,\,F)$, denote by $P_c$
the linearized Poincar\'{e} map of $c$. Then $P_c\in \Sp(2n-2)$ is symplectic.
For any $M\in \Sp(2k)$, we define the {\it elliptic height } $e(M)$
of $M$ to be the total algebraic multiplicity of all eigenvalues of
$M$ on the unit circle $\U=\{z\in\C|\; |z|=1\}$ in the complex plane
$\C$. Since $M$ is symplectic, $e(M)$ is even and $0\le e(M)\le 2k$.
A closed geodesic $c$ on the $n$-dimensional manifold $(M,F)$ is called {\it elliptic} if $e(P_c)=2(n-1)$, i.e., all the
eigenvalues of $P_c$ locate on $\U$; {\it hyperbolic} if $e(P_c)=0$, i.e., all the
eigenvalues of $P_c$ locate away from $\U$;
{\it non-degenerate} if $1$ is not an eigenvalue of $P_c$. A Finsler manifold $(M,\,F)$
is called {\it bumpy} if all the closed geodesics on it are non-degenerate.

There is a famous conjecture in Riemannian geometry which claims there exist infinitely many
closed geodesics on any compact Riemannian manifold. This conjecture
has been proved except for compact rank one symmetric
spaces. The results of Franks \cite{Fra} in 1992 and Bangert \cite{Ban}
in 1993 imply this conjecture is true for any Riemannian 2-sphere.
But once one move to the Finsler case, the conjecture
becomes false. It was quite surprising when Katok \cite{Kat}  in 1973 found
some irreversible Finsler metrics on spheres with only finitely
many closed geodesics and all closed geodesics are non-degenerate and elliptic (cf. \cite{Zil}).

Recently, index iterated theory of closed geodesics (cf. \cite{Bot} and \cite{Lon3}) has been applied to study the closed geodesic problem on Finsler manifolds. For example, Bangert and Long in \cite{BaL} show that there exist at least two closed geodesics on every $(S^2,F)$. After that, a great number of multiplicity and stability results has appeared (cf. \cite{DuL1}-\cite{DuL4}, \cite{DLW}, \cite{Lon4}, \cite{LoD}, \cite{LoW},  \cite{Rad3}, \cite{Wan1}-\cite{Wan4} and therein).

In \cite{Rad1}, Rademacher has firstly introduced the reversibility $\lambda=\lambda(M,F)$ of a compact Finsler manifold defined by
\bea \lambda=\max\{F(-X)\ |\ X\in TM,\ F(X)=1\}\ge 1.\nn\eea
Then Rademacher in \cite{Rad2} present some results about multiplicity and the length of closed geodesics and about their stability properties. For example, let $F$ be a Finsler metric on $S^{n}$ with reversibility $\lm$ and flag curvature $K$ satisfying $\left(\frac{\lm}{1+\lm}\right)^2<K\le 1$, then there exist at least $n/2-1$ closed geodesics with length $<2n\pi$. If $\frac{9\lm^2}{4(1+\lm)^2}<K\le 1$ and $\lm<2$, then there exists a closed geodesic of elliptic-parabolic, i.e., its linearized Poincar\'{e} map splits into $2$-dimensional rotations and a part whose eigenvalues are $\pm 1$. Some similar results in the Riemannian case (in this case, $\lambda=1$) have been obtained in \cite{BTZ1} and \cite{BTZ2}.

Recently, Wang in \cite{Wan1} proved that for every Finsler $n$-dimensional sphere $S^n$ with reversibility $\lm$ and flag curvature $K$ satisfying $\left(\frac{\lm}{1+\lm}\right)^2<K\le 1$, either there exist infinitely many prime closed geodesics or there exists one elliptic closed geodesics whose linearized Poincar\'{e} map has at least one eigenvalue which is of the form $\exp(\pi i\mu)$ with an irrational $\mu$. Furthermore, assume that this metric $F$ is bumpy, in \cite{Wan2}, Wang shows that there exist at least $2[\frac{n+1}{2}]$ closed geodesics on $(S^n,F)$. Also in \cite{Wan2}, Wang shows that for every bumpy Finsler metric $F$ on $S^n$ satisfying $\frac{9\lm^2}{4(1+\lm)^2}<K\le 1$, there exist two prime elliptic closed geodesics provided the number of closed geodesics on $(S^n,F)$ is finite.

In this paper, we generalize the above corresponding main results in \cite{Wan1}-\cite{Wan2} and \cite{Rad2} to the following theorem.

\medskip

{\bf Theorem 1.1.} {\it For a Finsler metric $F$ on the $n$-dimensional sphere $S^{n}$ with reversibility $\lm$ and flag curvature $K$ satisfying $\left(\frac{\lm}{1+\lm}\right)^2<K\le 1$, either there exist infinitely many closed geodesics, or there exist at least two elliptic closed geodesics and each linearized Poincar\'{e} map has at least one eigenvalue of the form $e^{\sqrt{-1}\th}$ with $\th$ being an irrational multiple of $\pi$.}

\medskip

Notice that this result is not new for $n=2$. Indeed, Long and Wang \cite{LoW} proved that for every Finsler $2$-dimensional sphere $(S^2,F)$, there exist at least two irrationally elliptic prime closed geodesics provided the number of prime closed geodesics is finite.

Our proof of Theorem 1.1 in Section 3 contains mainly three ingredients: the common index jump theorem of \cite{LoZ}, Morse theory and some new symmetric information about index jump. In applications of Morse theory, we follow some ideas from \cite{Wan1}-\cite{Wan2}. So, for the sake of shortening the length of this paper, we omitted lists of some well-known results about critical modules of closed geodesics, the mean index identity and Betti numbers of the free loop space $\ol{\Lm}S^n$ which can be found in \cite{LoD} and \cite{Wan1}-\cite{Wan2}. In addition, we also follow some ideas from our recent preprint \cite{DuL4}.

In this paper, let $\N$, $\N_0$, $\Z$, $\Q$, $\R$, and $\C$ denote
the sets of natural integers, non-negative integers, integers,
rational numbers, real numbers, and complex numbers respectively.
We define the functions
\be \left\{\matrix{[a]=\max\{k\in\Z\,|\,k\le a\}, &
           E(a)=\min\{k\in\Z\,|\,k\ge a\} , \cr
    \varphi(a)=E(a)-[a], &\{a\}=a-[a].  \cr}\right. \lb{1.1}\ee
Especially, $\varphi(a)=0$ if $ a\in\Z\,$, and $\varphi(a)=1$ if $
a\notin\Z\,$.

\setcounter{figure}{0}
\setcounter{equation}{0}
\section{Index iteration theory of closed geodesics}

In \cite{Lon1} of 1999, Long has established the basic normal form
decomposition of symplectic matrices. Based on this result he has
further established the precise iteration formulae of indices of
symplectic paths in \cite{Lon2} of 2000. Note that this index iteration formulae works for Morse indices of iterated closed geodesics on spheres (cf. \cite{Liu} and Chap. 12 of \cite{Lon3}). Since every closed geodesic on a sphere must be orientable. Then by Theorem 1.1 of \cite{Liu}, the initial Morse index of a closed geodesic on a Finsler $S^n$ coincides with the index of a corresponding symplectic path.

As in \cite{Lon2}, denote by
\bea
N_1(\lm, b) &=& \left(\matrix{\lm & b\cr
                                0 & \lm\cr}\right), \qquad {\rm for\;}\lm=\pm 1, \; b\in\R, \lb{2.1}\\
D(\lm) &=& \left(\matrix{\lm & 0\cr
                      0 & \lm^{-1}\cr}\right), \qquad {\rm for\;}\lm\in\R\bs\{0, \pm 1\}, \lb{2.2}\\
R(\th) &=& \left(\matrix{\cos\th & -\sin\th \cr
                           \sin\th & \cos\th\cr}\right), \qquad {\rm for\;}\th\in (0,\pi)\cup (\pi,2\pi), \lb{2.3}\\
N_2(e^{\th\sqrt{-1}}, B) &=& \left(\matrix{ R(\th) & B \cr
                  0 & R(\th)\cr}\right), \qquad {\rm for\;}\th\in (0,\pi)\cup (\pi,2\pi)\;\; {\rm and}\; \nn\\
        && \qquad B=\left(\matrix{b_1 & b_2\cr
                                  b_3 & b_4\cr}\right)\; {\rm with}\; b_j\in\R, \;\;
                                         {\rm and}\;\; b_2\not= b_3. \lb{2.4}\eea
Here $N_2(e^{\th\sqrt{-1}}, B)$ is non-trivial if $(b_2-b_3)\sin\theta<0$, and trivial
if $(b_2-b_3)\sin\theta>0$.

As in \cite{Lon2}, the $\diamond$-sum (direct sum) of any two real matrices is defined by
$$ \left(\matrix{A_1 & B_1\cr C_1 & D_1\cr}\right)_{2i\times 2i}\diamond
      \left(\matrix{A_2 & B_2\cr C_2 & D_2\cr}\right)_{2j\times 2j}
=\left(\matrix{A_1 & 0 & B_1 & 0 \cr
                                   0 & A_2 & 0& B_2\cr
                                   C_1 & 0 & D_1 & 0 \cr
                                   0 & C_2 & 0 & D_2}\right). $$

For every $M\in\Sp(2n-2)$, the homotopy set $\Omega(M)$ of $M$ in $\Sp(2n-2)$ is defined by
$$ \Om(M)=\{N\in\Sp(2n-2)\,|\,\sg(N)\cap\U=\sg(M)\cap\U\equiv\Gamma\;\mbox{and}
                    \;\nu_{\om}(N)=\nu_{\om}(M),\ \forall\om\in\Gamma\}, $$
where $\sg(M)$ denotes the spectrum of $M$,
$\nu_{\om}(M)=\dim_{\C}\ker_{\C}(M-\om I)$ for $\om\in\U$. Next we write $\nu(M)=\nu_1(M)$ when $\om=1$.
The component $\Om^0(M)$ of $P$ in $\Sp(2n-2)$ is defined by
the path connected component of $\Om(M)$ containing $M$.

\medskip

{\bf Theorem 2.1.} (cf. Theorem 7.8 of \cite{Lon1}, Theorems 1.2 and 1.3 of \cite{Lon2}, cf. also
Theorem 1.8.10, Lemma 2.3.5 and Theorem 8.3.1 of \cite{Lon3}) {\it For every $P\in\Sp(2n-2)$, there
exists a continuous path $f\in\Om^0(P)$ such that $f(0)=P$ and
\bea f(1)
&=& N_1(1,1)^{\dm p_-}\,\dm\,I_{2p_0}\,\dm\,N_1(1,-1)^{\dm p_+}
  \dm\,N_1(-1,1)^{\dm q_-}\,\dm\,(-I_{2q_0})\,\dm\,N_1(-1,-1)^{\dm q_+} \nn\\
&&\dm\,N_2(e^{\aa_{1}\sqrt{-1}},A_{1})\,\dm\,\cdots\,\dm\,N_2(e^{\aa_{r_{\ast}}\sqrt{-1}},A_{r_{\ast}})
  \dm\,N_2(e^{\bb_{1}\sqrt{-1}},B_{1})\,\dm\,\cdots\,\dm\,N_2(e^{\bb_{r_{0}}\sqrt{-1}},B_{r_{0}})\nn\\
&&\dm\,R(\th_1)\,\dm\,\cdots\,\dm\,R(\th_{r'})\,\dm\,R(\th_{r'+1})\,\dm\,\cdots\,\dm\,R(\th_r)\dm\,H(2)^{\dm h},\lb{2.5}\eea
where $\frac{\th_{j}}{2\pi}\in\Q\cap(0,1)$ for $1\le j\le r'$ and
$\frac{\th_{j}}{2\pi}\not\in\Q\cap(0,1)$ for $r'+1\le j\le r$; $N_2(e^{\aa_{j}\sqrt{-1}},A_{j})$'s
are nontrivial and $N_2(e^{\bb_{j}\sqrt{-1}},B_{j})$'s are trivial, and non-negative integers
$p_-, p_0, p_+,q_-, q_0, q_+,r,r_\ast,r_0,h$ satisfy
\be p_- + p_0 + p_+ + q_- + q_0 + q_+ + r + 2r_{\ast} + 2r_0 + h = n-1. \lb{2.6}\ee

Let $\ga\in\P_{\tau}(2n-2)=\{\ga\in C([0,\tau],\Sp(2n-2))\,|\,\ga(0)=I\}$. Denote the basic normal form
decomposition of $P\equiv \ga(\tau)$ by (\ref{2.5}). Then we have
\bea i(\ga^m)
&=& m(i(\ga)+p_-+p_0-r ) + 2\sum_{j=1}^r{E}\left(\frac{m\th_j}{2\pi}\right) - r   \nn\\
&&  - p_- - p_0 - {{1+(-1)^m}\over 2}(q_0+q_+)
              + 2\sum_{j=1}^{r_{\ast}}\vf\left(\frac{m\aa_j}{2\pi}\right) - 2r_{\ast}, \lb{2.7}\\
\nu(\ga^m)
 &=& \nu(\ga) + {{1+(-1)^m}\over 2}(q_-+2q_0+q_+) + 2\vs(m,\ga(\tau)),    \lb{2.8}\eea
where we denote by }
\be \vs(m,\ga(\tau)) = r - \sum_{j=1}^r\vf(\frac{m\th_j}{2\pi})
             + r_{\ast} - \sum_{j=1}^{r_{\ast}}\vf(\frac{m\aa_j}{2\pi})
             + r_0 - \sum_{j=1}^{r_0}\vf(\frac{m\bb_j}{2\pi}).    \lb{2.9}\ee

\medskip

The following is the common index jump theorem of Long and Zhu in \cite{LoZ}.

\medskip

{\bf Theorem 2.2.} (cf. Theorems 4.1-4.3 of \cite{LoZ}) {\it   Let $\gamma_k, k = 1,\ldots,q$ be a finite collection of symplectic paths. Let $M_k = \gamma_k(\tau_k)\in Sp(2n-2)$. Suppose the mean index $\hat i(\gamma_k, 1)=\lim_{m\rightarrow +\infty}\frac{i(\ga_k,m)}{m} > 0$, for any $k = 1,\ldots,q$.
Then there exist infinitely many $(N,m_1,\ldots,m_q)\in\N^{q+1}$  such that
\bea \nu(\gamma_k, 2m_k -1)&=&  \nu(\gamma_k, 1),\lb{2.10}\\
\nu(\gamma_k, 2m_k +1)&=&   \nu(\gamma_k, 1),\lb{2.11}\\
i(\gamma_k, 2m_k -1)+\nu(\gamma_k, 2m_k -1)&=&
2N-\left(i(\gamma_k, 1)+2S^+_{M_k}(1)-\nu(\gamma_k, 1)\right),\lb{2.12}\\
i(\gamma_k, 2m_k+1)&=&2N+i(\gamma_k, 1),\lb{2.13}\\
i(\gamma_k, 2m_k)&\ge&2N-\frac{e(M_k)}{2},\lb{2.14}\\
i(\gamma_k, 2m_k)+\nu(\gamma_k, 2m_k)&\le&2N+\frac{e(M_k)}{2},\lb{2.15}
\eea
for every $k=1,\ldots,q$, where $S^+_{M_k}(1)$ is the splitting number of $M_k$ at $1$.
Moreover we have
\bea \min\left\{\left\{\frac{m_k\theta}{\pi}\right\},\,1-\left\{\frac{m_k\theta}{\pi}\right\}\right\}
<\delta,\lb{2.16}\eea
whenever $e^{\sqrt{-1}\theta}\in\sigma(M_k)$
and $\delta$ can be chosen as small as we want (cf. (4.43) of \cite{LoZ}). More precisely, by (4.10) and (4.40) in \cite{LoZ} , we have
\bea m_k=\left(\left[\frac{N}{M\hat i(\gamma_k, 1)}\right]+\chi_k\right)M,\quad 1\le k\le q,\lb{2.17}\eea
where $\chi_k=0$ or $1$ for $1\le k\le q$ and $\{\frac{m_k\theta}{\pi}\}=0$
whenever $e^{\sqrt{-1}\theta}\in\sigma(M_k)$ and $\frac{\theta}{\pi}\in\Q$
for some $1\le k\le q$.}

\setcounter{figure}{0}
\setcounter{equation}{0}
\section{Proof of Theorem 1.1}

In order to prove Theorem 1.1, we make the following assumption

{\bf (FCG)} {\it Suppose that there exist only finitely many prime closed geodesics $\{c_k\}_{k=1}^q$, $k=1,\cdots,q$ on $(S^n, F)$ satisfying
$\left(\frac{\lambda}{\lambda+1}\right)^2<K\le 1$.}

\medskip

If the flag curvature $K$ of $(S^n, F)$ satisfies
$\left(\frac{\lambda}{\lambda+1}\right)^2<K\le 1$,
then every non-constant closed geodesic $c$ must satisfy
\bea i(c)&\ge& n-1, \lb{3.1}\\
\hat i(c)&>& n-1, \lb{3.2}\eea
where (\ref{3.1}) follows from Theorem 3 and Lemma 3 of \cite{Rad1},
(\ref{3.2}) follows from Lemma 2 of \cite{Rad2}.
Thus it follows from Theorem 2.2 of \cite{LoZ} (or, Theorem 10.2.3 of \cite{Lon3}) that
\bea i(c^{m+1})-i(c^m)-\nu(c^m)\ge i(c)-\frac{e(P_c)}{2}\ge 0,\quad\forall m\in\N,\lb{3.3}\eea
where the last inequality holds by (\ref{3.1}) and the fact that $e(P_c)\le 2(n-1)$.

By (\ref{3.2}) it follows from Theorem 2.2 that there exist infinitely many $(q+1)$-tuples $(N, m_1, \cdots, m_q)\in\N^{q+1}$ such that for any $1\le k\le q$, there holds
\bea i(c_k^{2m_k -1})+\nu(c_k^{2m_k -1})&=&
2N-\left(i(c_k)+2S^+_{M_k}(1)-\nu(c_k)\right), \lb{3.4}\\
i(c_k^{2m_k})&\ge& 2N-\frac{e(P_{c_k})}{2},\lb{3.5}\\
i(c_k^{2m_k})+\nu(c_k^{2m_k})&\le& 2N+\frac{e(P_{c_k})}{2},\lb{3.6}\\
i(c_k^{2m_k+1})&=&2N+i(c_k).\lb{3.7}\eea
where, together with (\ref{3.1}), (\ref{3.3}) and the fact $e(P_{c_k})\le 2(n-1)$, $1\le k\le q$, yields
\bea i(c_k^{m})+\nu(c_k^m)&\le& i(c_k^{2m_k}), \qquad\forall\ 1\le m<2m_k,\lb{3.8}\\
i(c_k^{2m_k})+\nu(c_k^{2m_k})&\le& 2N+\frac{e(P_{c_k})}{2}\le 2N+(n-1),\lb{3.9}\\
2N+(n-1)&\le &i(c_k^{m}),\qquad\forall\ m>2m_k.\lb{3.10}\eea

\medskip

{\bf Claim 1:} {\it Under the condition (FCG), there exists at least one closed geodesic among $\{c_k\}_{1\le k\le q}$, without loss of generality, saying $c_1$, satisfying $i(c_1^{2m_1})+\nu(c_1^{2m_1})=2N+(n-1)$.}

\medskip

If $i(c_k^{2m_k})+\nu(c_k^{2m_k})\le 2N+(n-2)$ for any $1\le k\le q$, then by (\ref{3.8})-(\ref{3.10}) we have
\bea i(c_k^{m})+\nu(c_k^m)&\le& i(c_k^{2m_k}), \quad\forall\ 1\le m<2m_k, \lb{3.11}\\
i(c_k^{2m_k})+\nu(c_k^{2m_k})&\le&2N+(n-2),\lb{3.12}\\
2N+(n-1)&\le&i(c_k^m),\quad \forall \ m>2m_k,\lb{3.13}\eea
So (6.9)-(6.11) in \cite{Wan1} hold. By Claims 1 and 2 in \cite{Wan1}, we obtain
\bea \sum_{p=0}^{2N+(n-2)}(-1)^p M_p=2N B(n,1)= \left\{\matrix{
             -\frac{N n}{n-1}, & \quad {\it if}\;\;n\in 2\N, \cr
              \frac{N(n+1)}{n-1}, & \quad {\it if}\;\;n\in 2\N-1, \cr}\right.\lb{3.14}\eea
where the Morse-type number $M_p$ is defined by $M_p =\sum_{1\le k\le q,\; m\ge 1}\dim{\ol{C}}_p(E, c^m_k)$ for all $p\in\Z$ and the critical module (cf. \cite{Cha}) is defined by
\bea \overline{C}_p( E,c_k^m)
\equiv H_p\left((\Lm(c_k^m)\cup S^1\cdot c_k^m)/S^1, \Lm(c_k^m)/S^1; \Q\right). \lb{3.15}\eea

Then the Morse inequality, together with (\ref{3.14}), yields a contradiction with the condition (FCG). Here we omitted the details because the proof is exactly the same as those in pages 1582-1585 of \cite{Wan1}. This completes the proof of Claim 1.

\medskip

In order to obtain more information about the closed geodesic $c_1$, we need to know the precise expression of $i(c_k^{2m_k})$ as follows
\bea i(c_k^{2m_k})
&=& 2m_k(i(c_k)+S^+_{M_k}(1)-C(M_k)) \nn\\
  &&+ 2\sum_{\th\in(0,2\pi)}E\left(\frac{m_k\th}{\pi}\right)S^-_{M_k}(e^{\sqrt{-1}\th})-(S^+_{M_k}(1)+C(M_k))\nn\\
&=& 2I(k,m_k)-(S^+_{M_k}(1)+C(M_k)) \nn\\
&=& 2(N+\Delta_k) -(S^+_{M_k}(1)+C(M_k)) \nn\\
&=& 2N -S^+_{M_k}(1)-C(M_k)+2\Delta_k,\qquad 1\le k\le q. \lb{3.16}\eea
where the first equality follows from Theorem 2.1 of \cite{LoZ} (cf. Theorem 9.3.1 of \cite{Lon3}), the second and third equalities follow from (4.38), (4.42) and
(4.45) in pages 349-350 of \cite{LoZ}, and denoted by
\bea \Delta_k = \sum_{0<\{m_k\th/\pi\}<\delta}S^-_{M_k}(e^{\sqrt{-1}\th}),\qquad
C(M_k) = \sum_{\th\in(0,2\pi)}S^-_{M_k}(e^{\sqrt{-1}\th}). \lb{3.17}\eea
Note that $\{\frac{m_k\th}{\pi}\}=0$ when $\frac{\th}{\pi}\in\Q$ by Theorem 2.2, so $\Delta_k$ only counts the number of those eigenvalues $e^{\sqrt{-1}\th}\in\sigma(M_k)$ and $\th$ being an irrational multiple of $\pi$ in the decomposition (\ref{2.5}) of $P_{c_k}$ satisfying $0<\{m_k\th/\pi\}<\delta$.

Notice that $\om=e^{\aa\sqrt{-1}}$ and its conjugate $\ol{\om}$ are bouble eigenvalues of nontrivial normal form $M=N_2(e^{\aa\sqrt{-1}},A)$ by Lemma 1.9.2 of \cite{Lon3}, so exactly half of eigenvalues $\{\om,\ol{\om}\}$'s of all nontrivial $N_2(e^{\aa\sqrt{-1}},A)$'s with $\frac{\aa}{\pi}\not\in\Q$ satisfy $0<\{m_k\aa/\pi\}<\delta$ and another half of ones satisfy $1-\delta<\{m_k\aa/\pi\}<1$. In addition, by List 9.1.12 of \cite{Lon3}, it yields $(S^-_M(\ol{\om}),S^-_M(\om))=(S^+_M(\om),S^-_M(\om))=(1,1)$. So by the definition (\ref{3.17}) of $\Delta_k$, there holds
\bea \Delta_k\le r_k-r_k'+r_{k_\ast}-r_{k_\ast}',\lb{3.18}\eea
where and below $r_k$ and $r_{k_\ast}$ denote the numbers of normal forms $R(\th)$ and nontrivial normal forms $N_2(e^{\aa\sqrt{-1}},A)$ in the decomposition (\ref{2.5}) of $P_{c_k}$, respectively, and $r_k'$ and $r_{k_\ast}'$ denote the numbers of rational normal forms $R(\th)$ and rational nontrivial normal forms $N_2(e^{\aa\sqrt{-1}},A)$, respectively.

\medskip

{\bf Claim 2:} {\it The closed geodesic $c_1$ found in Claim 1 is elliptic one whose linearized Poincar\'{e} map $P_{c_1}$ has at least one eigenvalue which is of the form $\exp(i\mu_1\pi)$ with an irrational $\mu_1$. Moreover, there holds
\bea r_1-r_1'=\Delta_1,\qquad r_1-r_1'\ge 1.\lb{3.19}\eea}
\indent Firstly by (\ref{3.9}) and the fact $i(c_1^{2m_1})+\nu(c_1^{2m_1})=2N+(n-1)$ from Claim 1, it yields $e(P_{c_1})=2(n-1)$, i.e., the closed geodesic $c_1$ is elliptic.

By (\ref{2.8}) and (\ref{2.9}) of Theorem 2.1 and the fact $\nu(c_k)=p_{k_-}+2p_{k_0}+p_{k_+}$ we have
\bea \nu(c_k^{2m_k})&=&\nu(c_k) + {{1+(-1)^{2m_k}}\over 2}(q_{k_-}+2q_{k_0}+q_{k_+}) \nn\\
     &&+ 2r_k - 2\sum_{j=1}^{r_k}\vf(\frac{m_k\th_j}{\pi})
             + 2r_{k_\ast} - 2\sum_{j=1}^{r_{k_\ast}}\vf(\frac{m_k\aa_j}{\pi})
             + 2r_{k_0} - 2\sum_{j=1}^{r_{k_0}}\vf(\frac{m_k\bb_j}{\pi})\nn\\
      &=& p_{k_-}+2p_{k_0}+p_{k_+}+q_{k_-}+2q_{k_0}+q_{k_+}+2r_k'+2r_{k_\ast}'+2r_{k_0}',\quad\forall\ 1\le k\le q,\lb{3.20}\eea
where and below $p_{k_-}, p_{k_0}, p_{k_+}, q_{k_-}, q_{k_0}, q_{k_+}, r_k, r_{k_\ast}, r_{k_0}, h_k$ denote the corresponding integers as those in the decomposition (\ref{2.5}) of $P_{c_k}$. And $r_{k_\ast}'$ and $r_{k_0}'$ denote the number of rational normal forms $N_2(e^{\aa\sqrt{-1}},A)$ and $N_2(e^{\bb\sqrt{-1}},B)$ in the decomposition (\ref{2.5}) of $P_{c_k}$, respectively.

Now by (\ref{3.16}), (\ref{3.20}) and the equality $i(c_1^{2m_1})+\nu(c_1^{2m_1})=2N+(n-1)$, it yields
\bea 2N+(n-1)&=&i(c_1^{2m_1})+\nu(c_1^{2m_1})\nn\\
    &=&2N-S^+_{M_1}(1)-C(M_1)+2\Delta_1\nn\\
      &&+p_{1_-}+2p_{1_0}+p_{1_+}+q_{1_-}+2q_{1_0}+q_{1_+}+2r_1'+2r_{1_\ast}'+2r_{1_0}'\nn\\
    &=& 2N-(p_{1_-}+p_{1_0})-(q_{1_0}+q_{1_+}+r_1+2r_{1_\ast})+2\Delta_1\nn\\
      &&+p_{1_-}+2p_{1_0}+p_{1_+}+q_{1_-}+2q_{1_0}+q_{1_+}+2r_1'+2r_{1_\ast}'+2r_{1_0}'\nn\\
    &=&2N+p_{1_0}+p_{1_+}+q_{1_-}+q_{1_0}+2r_{1_0}'-2(r_{1_\ast}-r_{1_\ast}')+2r_1'-r_1+2\Delta_1,\lb{3.21}\eea
where the third equality follows from $S^+_{M_1}(1)=p_{1_-}+p_{1_0}$ and $C(M_1)=q_{1_0}+q_{1_+}+r_1+2r_{1_\ast}$ by some results about splitting numbers in Lemma 9.1.5 and List 9.1.12 of \cite{Lon3}.

Then we rewrite (\ref{3.21}) and (\ref{2.6}) as follows
\bea p_{1_0}+p_{1_+}+q_{1_-}+q_{1_0}+2r_{1_0}'-2(r_{1_\ast}-r_{1_\ast}')+2r_1'-r_1+2\Delta_1&=&n-1,\lb{3.22}\\
p_{1_-} + p_{1_0} + p_{1_+} + q_{1_-} + q_{1_0} + q_{1_+} + r_1 + 2r_{1_\ast} + 2r_{1_0} + h_1 &=&n-1.\lb{3.23}\eea

So comparing with these two equalities, it yields
\bea p_{1_-} + q_{1_+} + 2(r_{1_\ast}-r_{1_\ast}'+r_1-r_1'-\Delta_1) + 2r_{1_\ast} + 2(r_{1_0}-r_{1_0}') + h_1=0.\lb{3.24}\eea

Note that every term in (\ref{3.24}) is a non-negative integer by (\ref{3.18}), so we obtain
\bea p_{1_-} = q_{1_+} = r_1-r_1'-\Delta_1 = r_{1_\ast} = r_{1_0}-r_{1_0}' = h_1=0.\lb{3.25}\eea

Firstly, if $\dim{\ol{C}}_{2N+(n-1)}(E, c_k^{2m_k})=0$, $1\le k\le q$, then the iterate $c_k^{2m_k}$ has
no contribution to the Morse type number $M_{2N+(n-1)}$. Then (\ref{3.14}) holds and the same argument in Claim 1 gives a contradiction.

Thus there exist at least $c_{k_0}$ such that $\dim{\ol{C}}_{2N+(n-1)}(E, c_{k_0}^{2m_{k_0}})\neq 0$. By the shifting property of critical modules (cf. Proposition 2.1 of \cite{LoD}) and (i) of Proposition 2.3 of \cite{LoD}, it yields $i(c_{k_0}^{2m_{k_0}})\le 2N+(n-1)\le i(c_{k_0}^{2m_{k_0}})+\nu(c_{k_0}^{2m_{k_0}})$ and $k_{\nu(c_{k_0}^{2m_{k_0}})}^{\bb(c_{k_0}^{2m_{k_0}})}\neq0$ (cf. Definition 2.2 of \cite{LoD} for notations $k_\ast^{\pm}(c)$'s). Together with (\ref{3.9}), this implies $i(c_{k_0}^{2m_{k_0}})+\nu(c_{k_0}^{2m_{k_0}})=2N+(n-1)$. Without loss of generality, we assume that $c_{k_0}$ is just $c_1$ since their ellipticity and index properties including (\ref{3.25}) are the same by above arguments.

In order to prove (\ref{3.19}), we assume that $r_1-r_1'=0$, i.e., all basic normal forms $R(\th)$'s in (\ref{2.5}) are rational, where, in this case, $c_k$ is called a rational closed geodesic in \cite{LoD}.  Note that $k_{\nu(c_1^{2m_1})}^{\bb(c_1^{2m_1})}\not\neq0$ by above arguments, let $g=c^{T(c_1)}$ (cf. Lemma 5.2 of \cite{Wan1} for the definition of $T(c_1)$). And an easy computation (note that $r_1=r_1'$ and (\ref{3.25}) is crucial in this proof, cf. Theorem 4.1 and its proof in \cite{LoD}) gives
\bea i(g^m)+\nu(g^m)=m(i(g)+\nu(g)-(n-1))+(n-1),\qquad\forall\ m\ge 1.\eea
Then the Hingston's result (cf. Proposition 1 of \cite{Hin1}, \cite{Hin2} and Theorem 4.2 of \cite{Wan1}) can applied to obtain the existence of infinitely many closed geodesics. See the details (6.37)-(6.40) in page 1857 of \cite{Wan1} (see Theorem 4.1 and its proof of \cite{LoD} as well) for this proof. Thus this also contradicts to the assumption (FCG).

So there holds $r_1-r_1'\ge 1$, i.e., there exist at least one eigenvalue of $P_{c_1}$ which is of the form $\exp(i\mu_1\pi)$ with an irrational $\mu_1$. This completes the proof of Claim 2.

\medskip

{\bf Claim 3:} {\it Under the condition (FCG), there exists another elliptic closed geodesic among $\{c_k\}_{2\le k\le q}$ and its linearized Poincar\'{e} map $P_{c_2}$ has at least one eigenvalue which is of the form $\exp(i\mu_2\pi)$ with an irrational $\mu_2$.}

\medskip

It follows from Theorem 2.2 that there exist also infinitely many $(q+1)$-tuples $(N', m_1', \cdots, m_q')\in\N^{q+1}$ such that for any $1\le k\le q$, together with (\ref{3.16}), there holds
\bea && i(c_k^{m})+\nu(c_k^m)\le i(c_k^{2m_k'}),\qquad \forall\ 1\le m<2m_k', \lb{3.27}\\
&&i(c_k^{2m_k'})+\nu(c_k^{2m_k'})\le 2N'+\frac{e(P_{c_k})}{2}\le 2N'+(n-1),\lb{3.28}\\
&& 2N'+(n-1)\le i(c_k^{m}),\qquad \forall\ m> 2m_k',\lb{3.29}\\
&& i(c_k^{2m_k'})=2N'-S^+_{M_k}(1)-C(M_k)+2\Delta_k',\lb{3.30}\eea
where, furthermore, $\Delta_k$ in (\ref{3.16}) and $\Delta_k'$ in (\ref{3.30}) satisfy the following relationship
\bea \Delta_k' + \Delta_k = r_k-r_k'+2(r_{k_\ast}-r_{k_\ast}'),\qquad\forall\ 1\le k\le q,\lb{3.31}\eea
where the equality follows from (especially the term (c) of) Theorem 4.2 of \cite{LoZ} or Lemma 3.18 and Corollary
3.19 of \cite{DuL3}. Then by (\ref{3.25}) it yields
\bea \Delta_1' + \Delta_1 =r_1-r_1'.\lb{3.32}\eea

At first, similar to the proof of (\ref{3.21}), by (\ref{2.5}), (\ref{3.19})-(\ref{3.20}), (\ref{3.25}) and (\ref{3.30})-(\ref{3.32}) we obtain
\bea i(c_1^{2m_1'})+\nu(c_1^{2m_1'})&=&2N'-S^+_{M_1}(1)-C(M_1)+2\Delta_1'\nn\\
   &&+p_{1_-}+2p_{1_0}+p_{1_+}+q_{1_-}+2q_{1_0}+q_{1_+}+2r_1'+2r_{1_\ast}'+2r_{1_0}'\nn\\
   &=& 2N'-(p_{1_-}+p_{1_0})-(q_{1_0}+q_{1_+}+r_1+2r_{1_\ast})+2(r_1-r_1')-2\Delta_1\nn\\
      &&+p_{1_-}+2p_{1_0}+p_{1_+}+q_{1_-}+2q_{1_0}+q_{1_+}+2r_1'+2r_{1_\ast}'+2r_{1_0}'\nn\\
    &=& 2N'-(p_{1_-}+p_{1_0})-(q_{1_0}+q_{1_+}+r_1+2r_{1_\ast})\nn\\
      &&+p_{1_-}+2p_{1_0}+p_{1_+}+q_{1_-}+2q_{1_0}+q_{1_+}+2r_1'+2r_{1_\ast}'+2r_{1_0}'\nn\\
    &=&2N'+p_{1_0}+p_{1_+}+q_{1_-}+q_{1_0}+2r_{1_0}'-2(r_{1_\ast}-r_{1_\ast}')+2r_1'-r_1\nn\\
    &=&2N'+p_{1_0}+p_{1_+}+q_{1_-}+q_{1_0}+2r_{1_0}'+r_1'-(r_1-r_1')\nn\\
    &\le & 2N'+(n-2),\lb{3.33}\eea
where we used (\ref{3.19}) and (\ref{3.25}) in the third and fifth equalities, respectively, and use (\ref{3.19}) and the fact 
$p_{1_0}+p_{1_+}+q_{1_-}+q_{1_0}+2r_{1_0}'+r_1'\le n-1$ in (\ref{2.5}) in the last inequality.

Assume that there holds $i(c_k^{2m_k'})+\nu(c_k^{2m_k'})<2N'+(n-1)$ for any $2\le k\le q$, then by (\ref{3.27})-(\ref{3.29}) and (\ref{3.33}) we get
\bea && i(c_k^{m})+\nu(c_k^m)\le i(c_k^{2m_k'}),\qquad \forall\ 1\le m<2m_k', \lb{3.34}\\
&&i(c_k^{2m_k'})+\nu(c_k^{2m_k'})\le 2N'+(n-2),\lb{3.35}\\
&& 2N'+(n-1)\le i(c_k^{m}),\qquad \forall\ m> 2m_k'.\lb{3.36}\eea

Thus as in Claim 1, we have
\bea \sum_{p=0}^{2N'+(n-2)}(-1)^p M_p=2N' B(n,1).\eea
Then the Morse inequality gives a contradiction.

Consequently there exist at least one closed geodesic among $\{c_k\}_{2\le k\le q}$, saying $c_2$, satisfying $i(c_2^{2m_2'})+\nu(c_2^{2m_2'})=2N'+(n-1)$.
Then the same proof as those in Claim 2 implies that $c_2$ is another elliptic closed geodesic whose linearized Poincar\'{e} map $P_{c_2}$ has at least one eigenvalue which is of the form $\exp(i\mu_2\pi)$ with an irrational $\mu_2$.

Now Claims 1-3 complete the proof of Theorem 1.1. \hfill\hb

\medskip

{\bf Acknowledgements.} The author would like to thank sincerely Professor Yiming Long for his encouragement, valuable comments, and helpful discussions with him. The author sincerely thank Professor Anatole Katok for his support and hospitality to him on his visit to the Department of Mathematics of Pennsylvania State University during August 2012 to August 2013. The author would like to thank Dr. Hui Liu reminding him that a missed term should be added to the right side of (\ref{3.18}).

\bibliographystyle{abbrv}

\end{document}